\documentclass[12pt]{article}
\usepackage{amsfonts}
\usepackage{mathrsfs}
\usepackage{amsmath,amssymb}
\baselineskip 5.5pt \lineskip 5.5pt \numberwithin{equation}{section}

\def\ma{\mathbb}
\def\BE {\begin{eqnarray}}
\def\EE {\end{eqnarray}}
\def\BC {\begin{eqnarray*}}
\def\EC {\end{eqnarray*}}
\def\OPLUS#1{\raisebox{-5pt}{\mbox{$\begin{array}{c}
\oplus\\[-5pt]\scriptstyle#1
\end{array}$}}}

\def\vs{\vspace*}\def\cl{\centerline}

\def\QED{\hfill$\Box$\par}
\def\a{\alpha}

\def\b{\beta}

\def\D{\Delta}

\def\Om{\Omega}

\def\LL{\mathcal {L}}

\def\cl{\centerline}

\def\ol{\overline}

\def\D{\Delta}

\def\bs{\backslash}

\def\vs{\vspace*}

\numberwithin{equation}{section}

\def\C{\mathbb{C}}
\def\Na{\nabla}

\def\Z{\mathbb{Z}}
\def\adddot{$\!\!\!${\bf.}\ \ }

\newtheorem{theo}{Theorem}[section]

\newtheorem{lemm}[theo]{Lemma}

\newtheorem{rema}[theo]{Remark}
\newtheorem{defi}[theo]{Definition}
\begin{document}
\cl{{\bf\large Classification of  modules of the intermediate
series}} \cl{\bf\large {over the twisted $N=2$ superconformal
algebra}\footnote {Supported by NSF grants 10671027, 10825101 of
China, ``One Hundred Talents Program'' from University of Science
and Technology of China.\\[2pt] \indent Corresponding E-mail:
sd\_junbo@163.com}} \vs{8pt}

\cl{Junbo Li$^{\dag,\ddag)}$, Yucai Su$^{\dag)}$, Linsheng
Zhu$^{\ddag)}$}

\cl{\small
$^{\dag)}$Department\,of\,Mathematics,\,University\,of\,Science\,and\,Technology\,of\,China,\,Hefei\,230026,\,China}

\cl{\small $^{\ddag)}$Department of Mathematics, Changshu Institute
of Technology, Changshu 215500, China}

\cl{\small E-mail: sd\_junbo@163.com, ycsu@ustc.edu.cn,
lszhu@cslg.edu.cn} \vs{8pt}

\noindent{\bf{Abstract.}} In this paper, a classification of modules
of the intermediate series over the twisted $N=2$ superconformal
algebra is obtained.

 \noindent{{\bf Key words:}
$N=2$ superconformal algebras, modules of the intermediate series,
Harish-Chandra modules.}

\noindent{\it{MR(2000) Subject Classification}: 17B10, 17B65,
17B68.}\vs{10pt}

\cl{\bf\S1. \
Introduction}\setcounter{section}{1}\setcounter{theo}{0}
\setcounter{equation}{0}

The $N=2$ superconformal algebras, constructed independently by Kac
{\cite{K1}} and by Ademollo et al.~{\cite{A1}}, whose structure and
representation theories were investigated in a series of papers
(e.g., \cite{DB,FJS,KL,KE}, etc), have attracted more and more
attentions in recent years. By choosing the standard Virasoro
generators, one observes that the superconformal algebras fall into
the following four types: (i) the Neveu-Schwarz algebra resulting
from antiperiodic boundary conditions, (ii) the Ramond algebra
corresponding to periodic boundary conditions, (iii) the topological
$N=2$ algebra  being as the symmetry algebra of topological
conformal field theory, (iv) the twisted superconformal algebra
possessing mixed boundary conditions. It is well known that the
first three superconformal algebras are isomorphic to each other.

The {\it twisted $N=2$ superconformal algebra} $\LL$ consists of the
Virasoro algebra generators $L_m,\,m\in\Z$ (corresponding to the
stress-energy tensor), the Heisenberg algebra generators $T_r$,
$r\in\frac{1}{2}+\Z$ (corresponding to the $U(1)$ current), the
fermionic generators $G_p,\,p\in\frac{1}{2}\Z$ (which are the modes
of the two spin-$\frac{3}{2}$ fermionic fields), with the following
commutation relations (where $c$ is a central element; see, e.g.,
\cite{DB})
\begin{eqnarray}\label{LieB}\begin{array}{lllll}
&& [L_m,
L_n]=(m-n)L_{n+m}+\frac{m^3-m}{12}\delta_{m+n,0}c,\vs{6pt}\\
&&
 [L_m,T_r]=-rT_{r+m},\ \ \ \ \ \ \ \ \ \ \ \ \ \ \ \ [T_r, T_s]=\frac{r}{3}\delta_{r+s,0}c,\vs{6pt}\\
&&
[L_m,G_p]=(\frac{m}{2}-p)G_{p+m},\ \ \ \ \ \ \ \ [T_r,G_p]=G_{p+r},\vs{6pt}\\
&&[G_p,G_q]=\left\{\begin{array}{ll}(-1)^{2p}\big(2L_{p+q}+\frac13(p^2-\frac14)\delta_{p+q,0}c\big)&\mbox{if
\
}p+q\in\Z,\\[6pt]
(-1)^{2p+1}(p-q)T_{p+q}&\mbox{if \
}p+q\in\frac12+\Z,\end{array}\right.
\end{array}\end{eqnarray}
for $m,n\in\Z,\,r,s\in\frac{1}{2}+\Z,\,p,q\in\frac{1}{2}\Z$.

Obviously, $\LL$ is ${\mathbb{Z}}_2$-graded: ${\cal L}={\cal
L}_{\overline{0}}\oplus{\cal L}_{\overline{1}},$ with
\begin{eqnarray}\label{0426a002}
{\cal L}_{\overline{0}}=\mbox{span}_{\ma C}\{
L_m,\,T_r,\,c\,|\,m\in\Z,\,r\in\frac{1}{2}+\Z\},\ \ \ \ {\cal
L}_{\overline{1}}=\mbox{span}_{\ma C}\{G_p\,|\,p\in\frac{1}{2}\Z\},
\end{eqnarray}
and the Cartan subalgebra  ${\cal H}={\mathbb C}L_0+{\mathbb C}c$.
One can easily see that $\LL$ contains  the well-known Virasoro
algebra $\mathcal {V}ir=$ span$_{\mathbb C} \{
L_m,c\,|\,m\in{\mathbb Z}\}$, the super-Virasoro algebras $\mathcal
{NS}$ (the $N=1$ Neveu-Schwarz algebra \cite{NS}) spanned by
$\{L_n,G_{r},c\,|\,n\in\Z,r\in\frac{1}{2}+\Z\}$, and $\mathcal {R}$
(the $N=1$ Ramond algebra \cite{R}) spanned by
$\{L_n,G_{n},c\,|\,n\in\Z\}$.

An ${\cal L}$-module $V=V_{\ol0}\oplus V_{\ol1}$ is called a {\it
Harish-Chandra module} if $V_\a$ is a direct sum of its finite
dimensional weight spaces $V_\a^\lambda=\{v\in V_\a\,|\,x
v=\lambda(x)v,\,x\in{\cal H}\}$ for $\a\in\Z_2$ and all
$\lambda\in{\cal H}^*$ (the dual space of ${\cal H}$). Similar to
the case of Virasoro algebra, we can define the module of the
intermediate series over $\LL$.
\begin{defi}\adddot
\rm A $\Z_2$-graded module $V=V_{\ol0}\oplus V_{\ol1}$ over ${\cal
L}$ is called a {\it module of the intermediate series} if $V$ is a
Harish-Chandra module such that ${\rm dim}V_\a^\lambda\leqslant1$
for $\lambda\in{\cal H}^*,\,\a\in\Z_2$.
\end{defi}

Although much progress on representation theory over the $N=2$
superconformal algebras has been made in recent years and many
aspects arise that are new to representations of chiral algebras,
there are still many questions  (e.g., structure of Verma modules,
irreducible modules, classification of Harish-Chandra modules,
etc.). In order to better understand the representations of the Lie
algebra $\cal L$, it is highly desirable to give a classification of
modules of the intermediate series over ${\cal L}$. This is the task
of the present paper.

The similar problem has been considered over the Ramond sector in
\cite{FJS}. However, due to the facts that the twisted sector is not
isomorphic to the Ramond sector and the number of fermionic currents
$N$ increases, the modules over the twisted sector are interesting
and complicated. It is well known that the twisted $N=2$ algebra
behaves very differently from the other $N=2$ algebras with respect
to the adapted ordering method and thus gives more new and different
insights in superconformal representation theory. In \cite{DB}, the
authors created the basis for the study of the embedding structure
of the twisted $N=2$ highest weight representations and deduced the
singular dimensions as well as the multiplication rules for singular
vector operators, which both are crucial for the analysis of the
embedding structure. In \cite{IK}, the structures of both Verma
modules and Fock modules over the twisted sector of the $N=2$
superconformal algebras were investigated. Because of lack of coset
construction of the twisted sector, the methods used in other
sectors are not applicable in the twisted case.

Now let us formulate the main result below. First we define two
families $A_{a,b},\,B_{a,b},\,a,b\in\C$ of indecomposable modules of
the intermediate series over $\cal L$ as follows: They all have
basis $\{x_n,y_n\,|\,n\in\frac12\Z\}$ such that (in all cases, $c$
acts trivially, and
$n\in\Z,\,r\in\frac{1}{2}+\Z,\,q,k\in\frac{1}{2}\Z$)
\begin{eqnarray}\label{modstr1}\begin{array}{llllll}
A_{a,b}: &L_nx_k=(a-k+bn)x_{n+k},&
L_ny_k=\big(a-k+(b+\frac{1}{2})n\big)y_{n+k},
\vs{4pt}\\
&T_rx_k=-2(b+1)x_{k+r},&
T_ry_k=-(2b+1)y_{k+r},\vs{4pt}\\
&G_qx_k=y_{q+k},& G_qy_k=(-1)^{2q}(a-k+2b q+q)x_{q+k},
\end{array}\end{eqnarray}
and
\begin{eqnarray}\label{modstr2}\begin{array}{llllll}
B_{a,b}:&L_nx_k=(a-k+bn)x_{n+k},\ \ \
L_ny_k=\left\{\begin{array}{ll}
\big(a-k+(b-\frac{1}{2})n\big)y_{n+k}&{\rm if}\ \,k\in\Z,\vs{4pt}\\
\big(a-k+(b+\frac{1}{2})n\big)y_{n+k}&{\rm if}\
\,k\in\frac{1}{2}+\Z,
\end{array}\right.\vs{4pt}\\
&T_rx_k=x_{k+r},\ \ \ \ \ \ \ \ \ \ \ \ \ \ \ \ \ \ \ \
T_ry_k=\left\{\begin{array}{ll}
0&{\rm if}\ \,k\in\Z,\vs{4pt}\\
(2b+1)\,r\,y_{k+r},&{\rm if}\ \,k\in\frac{1}{2}+\Z,
\end{array}\right.\vs{4pt}\\
&G_qx_k=\left\{\begin{array}{cl}
(-1)^{2q}(a-k+2bq)y_{k+q}&{\rm if}\ \,q+k\in\Z,\vs{4pt}\\
(-1)^{2q+1}y_{k+q}&{\rm if}\ \,q+k\in\frac{1}{2}+\Z,
\end{array}\right.\vs{4pt}\\
&G_qy_k=\left\{\begin{array}{cl}
x_{k+q}&{\rm if}\ \,k\in\Z,\vs{5pt}\\
-(a-k+2bq+q)x_{k+q}&{\rm if}\ \,k\in\frac{1}{2}+\Z,
\end{array}\right.
\end{array}\end{eqnarray}
The modules $A_{0,-1},\,A_{0,-\frac12}$ have deformations denoted by
$A_1(\a),\,A_2(\a)$, $\a\in\C$, defined by
\begin{eqnarray}\label{modstr3}
&\!\!\!\!\!\!\!\!\!\!\!\!\!\!\!\!\!\!\!\!\!&
\begin{array}{llllll}
A_1(\a):\!\!\!\!&&L_nx_k=\left\{\begin{array}{cl}\!\!-(k+n)x_{n+k}&{\rm
if}\
\,k\in\frac12\Z^*,\vs{4pt}\\
\!\!-n(n+\a)x_{n+k}&{\rm if}\ \,k=0,\end{array}\right.&
L_ny_k=-(k+\frac{n}{2})y_{n+k},\vs{6pt}\\
&&T_rx_k=\left\{\begin{array}{cl}
\!\!0&{\rm if}\ \,k\in\frac12\Z^*,\vs{6pt}\\
\!\!-2rx_{r+k}&{\rm if}\ \,k=0,
\end{array}\right.&
T_ry_k=y_{k+r},\vs{6pt}\\
&&G_qx_k=\left\{\begin{array}{cl}\!\!y_{q+k}&{\rm if}\
\,k\in\frac12\Z^*,\vs{4pt}\\
\!\!(2p+\alpha)y_{q+k}&{\rm if}\ \,k=0,
\end{array}\right.&
G_qy_k=(-1)^{2q+1}(k+q)x_{q+k},
\end{array}
\\%\end{eqnarray}\vs{-22pt}\begin{eqnarray}
\label{modstr4} &\!\!\!\!\!\!\!\!\!\!\!\!\!\!\!\!\!\!\!\!\!&
\begin{array}{llllll}
A_2(\a):\!\!\!\!&&L_nx_k=-(k+\frac{n}{2})x_{n+k},&
L_ny_k=\left\{\begin{array}{cl} -ky_{n+k}&{\rm if}\
\,k\in\frac12\Z\bs\{-n\},\vs{4pt}\\
n(n+\a)y_{0}&{\rm if}\ \,k=-n,
\end{array}\right.\vs{6pt}\\
&&T_rx_k=-x_{k+r},& T_ry_k=\left\{\begin{array}{cl}
0&{\rm if}\ \,k\in\frac12\Z\bs\{-r\},\vs{6pt}\\
2ry_{0}&{\rm if}\ \,k=-r,
\end{array}\right.\vs{6pt}\\
&&G_qy_k=(-1)^{2q+1}kx_{q+k},& G_qx_k=\left\{\begin{array}{cl}
y_{q+k}&{\rm if}\
\,k\in\frac12\Z\bs\{-q\},\vs{4pt}\\
(2p+\alpha)y_{0}&{\rm if}\ \,k=-q.
\end{array}\right.
\end{array}\end{eqnarray}
Similarly, the modules $B_{0,-\frac12},\,B_{\frac12,-\frac12}$ have
deformations denoted by $B_1(\a),\,B_2(\a)$,
\begin{eqnarray}\label{modstr66}\begin{array}{llllll}
B_1(\a):\!\!\!\!&&L_nx_k=-(k+\frac{n}{2})x_{n+k},\ \
L_ny_k=\left\{\begin{array}{cl}
-(k+n)y_{n+k}&{\rm if}\ \,k\in\Z^*,\vs{4pt}\\
-n(n+\a)y_{n+k}&{\rm if}\ \,k=0,\vs{4pt}\\
-ky_{n+k}&{\rm if}\ \,k\in\frac12+\Z,
\end{array}\right.\vs{6pt}\\
&&T_rx_k=x_{k+r},\ \ \ \ \ \ \ \ \ \ \ \ \ \ \ \,
T_ry_k=\left\{\begin{array}{cl}
0&{\rm if}\ \,k\in\frac12\Z^*,\vs{4pt}\\
-2y_{r+k}&{\rm if}\ \,k=0,
\end{array}\right.\vs{4pt}\\
&&G_qx_k=\left\{\begin{array}{cl}
(-1)^{2q+1}(k+q)y_{k+q}&{\rm if}\ \,q+k\in\Z,\vs{4pt}\\
(-1)^{2q+1}y_{k+q}&{\rm if}\ \,q+k\in\frac{1}{2}+\Z,
\end{array}\right.\vs{6pt}\\
&&G_qy_k=\left\{\begin{array}{cl}
x_{q+k}&{\rm if}\ \,k\in\Z^*,\vs{4pt}\\
(2p+\alpha)x_{q+k}&{\rm if}\ \,k=0,\vs{4pt}\\
kx_{k+q}&{\rm if}\ \,k\in\frac{1}{2}+\Z,
\end{array}\right.
\end{array}\end{eqnarray}

\begin{eqnarray}\label{modstr68}\begin{array}{llllll}
B_2(\a):\!\!\!\!\!\!&& L_nx_k=(\frac{1}{2}-k-\frac{1}{2}n)x_{n+k},\,
L_ny_k\!=\!\left\{\begin{array}{cl}
\!\!(\frac{1}{2}-k-n)y_{n+k}&\!\!\!{\rm if}\ \,k\in\Z,\vs{6pt}\\
\!\!(\frac{1}{2}-k)y_{n+k}&\!\!\!{\rm if}\
\,k\in\frac{1}{2}+(\Z\bs\{-n\}),\vs{6pt}\\
\!\!n(n+\a)y_{n+k}&\!\!\!{\rm if}\ \,k=\frac{1}{2}-n,
\end{array}\right.\vs{9pt}\\
&&T_rx_k=x_{k+r},\ \ \ \ \ \ \ \ \ \ \ \ \ \ \ \ \ \
T_ry_k=\left\{\begin{array}{cl} 0&{\rm if}\
\,k\in\frac{1}{2}\Z\bs\{\frac{1}{2}-r\},\vs{6pt}\\
2y_{r+k}&{\rm if}\ \,k=\frac{1}{2}-r,
\end{array}\right.\vs{9pt}\\
&&G_qx_k=\left\{\begin{array}{cl}
(-1)^{2q}(\frac{1}{2}-k-q)y_{k+q}&{\rm if}\ \,q+k\in\Z,\vs{6pt}\\
(-1)^{2q+1}y_{k+q}&{\rm if}\
\,q+k\in\frac{1}{2}+(\Z\bs\{-q\}),\vs{6pt}\\
(-1)^{2q+1}(2q+\alpha)y_{k+q}&{\rm if}\ \,k=\frac{1}{2}-q,
\end{array}\right.\vs{9pt}\\
&&G_qy_k=\left\{\begin{array}{cl}
x_{k+q}&{\rm if}\ \,k\in\Z,\vs{6pt}\\
(k-\frac{1}{2})x_{k+q}&{\rm if}\ \,k\in\frac{1}{2}+\Z.
\end{array}\right.
\end{array}\end{eqnarray}

The main result can be formulated as the following theorem.
\begin{theo}\adddot\label{mainth}
Any indecomposable module of the intermediate series $V$ over the
twisted $N=2$ superconformal algebra is one of the modules
$A_{a,b}$, $B_{a,b}$, $A_1(\a )$, $A_2(\a )$, $B_1(\a )$, $B_2(\a
)$, or one of their quotient modules for some  $a, b,\a\in{\C}$.
\end{theo}

\vs{5pt}

\cl{\bf\S2. \ The indecomposable Harish-Chandra module
$A_{a,b}$}\setcounter{section}{2}\setcounter{theo}{0}
\setcounter{equation}{0}

\vs{4pt}

It is well known that a module of the intermediate series over the
super-Virasoro algebra $\mathcal {NS}$ is one of the three series of
the modules $\mathcal {A}_{a,b},\,\mathcal {A}(\a),\,\mathcal
{B}(\b)$ or their quotient modules for suitable $a,b,\a,\b\in\C$,
where $\mathcal {A}_{a,b},\,\mathcal {A}(\a)$ have basis
$\{x_k\,|\,k\in\Z\}\cup\{y_u\,|\,u\in\frac{1}{2}+\Z\}$
 and $\mathcal {B}(\b)$ has basis $\{x_u\,|\,u\in\frac{1}{2}+\Z\}\cup\{y_k\,|\,k\in\Z\}$ such
that $c$ acts trivially and  (see, e.g., \cite{S1,SZ,S3})
\begin{eqnarray*}
{\mathcal {A}}_{a,b}:\,\,L_nx_k\!\!\!&=&\!\!\!(a-k+bn)x_{n+k},\ \
L_ny_u=\big(a-u+n(b+\frac{1}{2})\big)y_{n+u},\\
\,\,G_px_k\!\!\!&=&\!\!\!y_{k+p},\ \ \ \ \ \ \ \ \ \ \ \ \ \,\ \ \ \
\
G_py_u=-\big(a-u+2p(b+\frac{1}{2})\big)x_{p+u},\\
\mathcal
{A}(\a):\,\,L_nx_k\!\!\!&=&\!\!\!-(k+n)x_{n+k},\,k\neq0,\,L_nx_0=-n(n+\a)x_{n},\,
L_ny_u=-(u+\frac{n}{2})y_{n+u},\\
\,\,G_px_k\!\!\!&=&\!\!\!y_{k+p},\ \ k\neq0,\ \
G_px_0=(2p+\a)y_{p},\ \
G_py_u=(u+p)x_{p+u},\\
\mathcal { B}(\b):\,\,L_nx_u\!\!\!&=&\!\!\!-(u+\frac{n}{2})x_{n+u},\
\ \ L_ny_k=-ky_{n+k},\,\ k\neq-n,\,\
L_ny_{-n}=n(n+\b)y_{0},\\
\,\,G_px_u\!\!\!&=&\!\!\!y_{p+u},\ \ u\neq-p,\ \
G_px_{-p}=(2p+\b)y_{0},\ \ G_py_k=kx_{p+k},
\end{eqnarray*}
for $n,k\in\Z,\,p,u\in\frac{1}{2}+\Z$.

Let $V=V_{\ol0}\oplus V_{\ol1}$ be any indecomposable ${\cal
L}$-module with ${\rm dim}V^{\lambda}_{\alpha}\leqslant1$ for all
$\lambda\in{\cal H}^*,\,\alpha\in\Z_2$, where
\begin{eqnarray*}
V^{\lambda}_{\alpha}=\{v\in V_{\alpha}\mid L_0\cdot
v=\lambda(L_0)v\}.
\end{eqnarray*}
It is easy to see that $c$ acts trivially on $V$ (see, e.g.,
\cite{IL, CA}). So we can omit $c$ in the following. Being
indecomposable, it is clear that there exists $a\in\C$ such that
\begin{eqnarray}
\label{a30}
V=\big(\OPLUS{k\in\frac{1}{2}\Z}V_{\ol0}^{a+k}\big)\oplus\big(\OPLUS{k\in\frac{1}{2}\Z}V_{\ol1}^{a+k}\big).
\end{eqnarray}
Obviously
$V'=\oplus_{k\in\Z}V_{\ol0}^{a+k}\oplus\oplus_{k\in\frac12+\Z}V_{\ol1}^{a+k}$
and
$V''=\oplus_{k\in\frac12+\Z}V_{\ol0}^{a+k}\oplus\oplus_{k\in\Z}V_{\ol1}^{a+k}$
are modules of the intermediate series over $\cal NS$. First we
assume that they are modules of type $\mathcal {A}_{a,b}$ (later on,
we shall consider all possible deformations). Thus we can choose a
basis $\{x_k,y_k\,|\,k\in\frac12\Z\}$ of $V$ such that $V_{\ol0}$ is
spanned by $\{x_i,y_{i+\frac{1}{2}}\,|\,i\in\Z\}$ and $V_{\ol1}$ is
spanned by $\{x_{i+\frac{1}{2}},y_i\,|\,i\in\Z\}$, with the
following relations, for some $a,\,b,\,b'\in{\ma C}$,
\begin{eqnarray}\label{modstr101}
\begin{array}{llllll}
&&L_mx_i=\left\{\begin{array}{ll}
(a-i+bm)x_{m+i}&{\rm if}\ \,i\in\Z,\vs{6pt}\\
(a-i+b'm)x_{m+i}&{\rm if}\ \,i\in\frac{1}{2}+\Z,
\end{array}\right.\vs{9pt}\\
&&L_my_j=\left\{\begin{array}{ll}
\big(a-j+m(b'+\frac{1}{2})\big)y_{m+j}&{\rm if}\ \,j\in\Z,\vs{6pt}\\
\big(a-j+m(b+\frac{1}{2})\big)y_{m+j}&{\rm if}\
\,j\in\frac{1}{2}+\Z,
\end{array}\right.\vs{6pt}\\
&&G_px_k=y_{k+p}\,,\\[6pt]
 &&G_py_i=\left\{\begin{array}{ll}
-\big(a-i+2p(b'+\frac{1}{2})\big)x_{p+i}&{\rm if}\ \,i\in\Z,\vs{6pt}\\
-\big(a-i+2p(b+\frac{1}{2})\big)x_{p+i}&{\rm if}\
\,i\in\frac{1}{2}+\Z,
\end{array}\right.
\end{array}\end{eqnarray}
where $m\in\Z,\,r,p\in\frac{1}{2}+\Z,\,i,j,k\in\frac{1}{2}\Z$. To
determine possible structures on $V$, we suppose
\begin{eqnarray}\label{modstr102}
\begin{array}{llllll}
&&T_rx_k=f_{r,k}x_{k+r},\ \,\,\ T_ry_k=f'_{r,k}y_{k+r},
\ \ \ %\vs{6pt}\\&&
G_nx_k=g_{n,k}y_{k+n},\ \,\ G_ny_k=g'_{n,k}x_{k+n},\,\
\end{array}\end{eqnarray}
for $n\in\Z,\,r\in\frac{1}{2}+\Z,\,k\in\frac{1}{2}\Z$ and some
$f_{r,k},\,f'_{r,k},\,g_{n,k},\,g'_{n,k}\in{\ma C}$. Denote such a
module defined in (\ref{modstr101}) and (\ref{modstr102}) by
$A_{a,b,b'}$ temporarily.
\begin{lemm}\adddot\label{lemma2.1}
{\rm(i)}\ \,If there is some $r_0\in\frac12+\Z,\,k_0\in\frac12\Z$
such that $f_{r_0,k_0}\neq0$, then for any fixed $n\in\Z$, there is
always infinitely many $m\in\Z$ satisfying $f_{r_0+n,k_0+m}\neq0$.
The analogous result also holds for $f'_{r,k},g_{n,k}$ and $g'_{n,k}$.\\
{\rm(ii)}\ \,For any $r_0\in\frac12+\Z,\,k_0\in\frac12\Z$, there are
infinitely many pairs $(n,m)\in\Z\times\Z$ with
$f_{r_0+n,k_0+m}\ne0$, or there are infinitely many pairs
$(n,m)\in\Z\times\Z$ with
$f'_{r_0+n,k_0+m}\ne0$. \\
{\rm(iii)}\ \,For any $n_0\in\Z,\,k_0\in\frac12\Z$, there are
infinitely many pairs $(n,m)\in\Z\times\Z$ with
$g_{n_0+n,k_0+m}\ne0.$ Similar result also holds for $g'_{n,k}$.
\end{lemm}
{\bf Proof.}\ \ We shall  prove  (i) for $g_{n,k}$. The rest can be
proved similarly. Suppose there is some $n_0,k_0\in\Z$ with
$g_{n_0,k_0}\neq0$. For $m,n,k\in\Z$, applying
$[L_m,G_n]=(\frac{m}{2}-n)G_{m+n}$ to $x_k$ gives
\begin{eqnarray}\label{0509a01}
&&g_{n,k}\big(a-n-k+m(b'+\frac{1}{2})\big)-(a-k+bm)g_{n,m+k}=(\frac{m}{2}-n)g_{m+n,k}.
\end{eqnarray}
According to (\ref{0509a01}), we only need to prove that for any
$n\in\Z$, there is some $k\in\Z$ such that $g_{n,k}\neq0$. Assume
conversely, say, $g_{1,k}=0$ for $k\in\Z$. By (\ref{0509a01}), we
obtain $g_{m,k}=0$ for $m\in\Z\bs\{3\},\,k\in\Z$. Taking
$m=4,\,n=-1$, one has $g_{3,k}=0$ for $k\in\Z$, which implies
$g_{m,k}=0$ for $m,k\in\Z$, a contradiction. The lemma follows.\QED
\begin{lemm}\adddot\label{lemma2.3} $b'=b$, where $b$ and $b'$ are defined in
$(\ref{modstr101}).$
\end{lemm}
{\bf Proof.}\ \ According to Lemmas \ref{lemma2.1}, we first suppose
\begin{eqnarray}\label{supp040701}
f_{r,k}\neq 0\ \ \mbox{for infinitely many pairs}\ \,
(r,k)\in(\frac{1}{2}+\Z)\times\Z.
\end{eqnarray}
For $m,\,n,\,k\in\Z,\,r\in\frac{1}{2}+\Z$, applying
$%\begin{eqnarray}\label{eqact2of1st}&&
\big[L_m,[L_n,T_r]\big]=-(n+r)[L_{m+n},T_r] $ %\end{eqnarray}
to
$x_k$ and comparing the coefficients of $x_{k+m+n+r}$, one has
\begin{eqnarray}
&&\!\!\!\!\!\!\!\!(a-k+bm)(a-k-m+bn)f_{r,k+m+n}-(a-k+bn)(a-k-n-r+b'm)f_{r,k+n}\nonumber\\
&&\!\!\!\!\!\!\!\!+(a-k-n-r+b'm)(a-k-r+b'n)f_{r,k}-(a-k+bm)(a-k-m-r+b'n)f_{r,k+m}\nonumber\\
&&\!\!\!\!\!\!\!\!=(n+r)(a-k+bm+bn)f_{r,k+m+n}-(n+r)(a-k-r+b'm+b'n)f_{r,k}.\label{eqLTTxI0}
\end{eqnarray}
Replacing $m,\,n,\,k$ by (i) $m,\,m,\,k-m$, (ii) $-m,\,-m,\,k+m$ and
(iii) $m,\,-m,\,k$ respectively, we obtain the following three
equations
\begin{eqnarray}
&&\!\!\!\!\!\!\!\!\!\!\!\!\!\!\!\!\!
\big((a-k+bm)(a-k+bm+m)-(m+r)(a-k+2bm+m)\big)f_{r,k+m}\nonumber\\
&&\!\!\!\!\!\!\!\!\!
-2(a-k+bk+m)(a-k-r+b'm)f_{r,k}+\big((m+r)(a-k-r\nonumber\\
&&\!\!\!\!\!\!\!\!\!
+2b'm+m)+(a-k-r+b'm)(a-k-r+m+b'm)\big)f_{r,k-m}=0,\label{eqLTTxI011}\\[4pt]
&&\!\!\!\!\!\!\!\!\!\!\!\!\!\!\!\!\!
\big((r-m)(a-k-m-r-2b'm)+(a-k-r-b'm)(a-k-r\nonumber\\
&&\!\!\!\!\!\!\!\!\!
-b'm-m)\big)f_{r,k+m}-2\big(a-k-(b+1)m\big)(a-k-r-b'm)f_{r,k}\nonumber\\
&&\!\!\!\!\!\!\!\!\!
+\big((m-r)(a-k-2bm-m)+(a-k-b m)(a-k-m-b
m)\big)f_{r,k-m}=0,\label{eqLTTxI02}\\[4pt]
&&\!\!\!\!\!\!\!\!\!\!\!\!\!\!\!\!\!
(k-a-b m)(a-k-r-m-b'm)f_{r,k+m}+\big((a-k-r-b'm)(a-k-r\nonumber\\
&&\!\!\!\!\!\!\!\!\!
+b'm+m)+(a-k-b m-m)(a-k+b m)+r(r-m)\big)f_{r,k}\nonumber\\
&&\!\!\!\!\!\!\!\!\!
+(k+bm-a)(a-k-r+b'm+m)f_{r,k-m}=0.\label{eqLTTxI03}
\end{eqnarray}
Regard (\ref{eqLTTxI011})--(\ref{eqLTTxI03}) as a system of 3 linear
equations on the variables $f_{r,k+m},\,f_{r,k},\,f_{r,k-m}$. Denote
by $\D_1(m)$ the determinant of coefficients (which is a polynomial
on $m$). By (\ref{supp040701}), we must have
\begin{equation} \label{Delta1=01}
\D_1(m)=0 \mbox{ \ for infinitely many thus for all \ }m\in\Z,
\end{equation}
where by a lengthy computation,
$$\D_1(m)=(b-b'-2)(b-b'-1)(b-b')(b+b'+1)(b^2+b+2bb'+3b'+{b'}^2)m^6.
$$
Thus (\ref{Delta1=01}) shows
\begin{eqnarray}\label{Delta12=0}
&&b'\in\{-b-1,b-2,b-1,b,\frac{1}{2}(-3-2b\pm \sqrt{9+8b})\}.
\end{eqnarray}
If (\ref{supp040701}) does not hold, then Lemma \ref{lemma2.1} shows
(\ref{supp040701}) must hold for $f'_{r,k}$, from this we obtain as
in (\ref{eqLTTxI0}),
\begin{eqnarray*}
&&\big(a-k+(b+\frac{1}{2})m\big)\big(a-k-m+(b+\frac{1}{2})n\big)f'_{r,k+m+n}\\
&&-\big(a-k+(b+\frac{1}{2})n\big)\big(a-k-n-r+(b'+\frac{1}{2})m\big)f'_{r,k+n}\nonumber\\
&&+\big(a-k-n-r+(b'+\frac{1}{2})m\big)\big(a-k-r+(b'+\frac{1}{2})n\big)f'_{r,k}\\
&&-\big(a-k+(b+\frac{1}{2})m\big)\big(a-k-m-r+(b'+\frac{1}{2})n\big)f'_{r,k+m}\nonumber\\
&&=(n+r)\big(a-k+(b+\frac{1}{2})m+(b+\frac{1}{2})n)f'_{r,k+m+n}\\
&&-(n+r)\big(a-k-r+(b'+\frac{1}{2})m+(b'+\frac{1}{2})n\big)f'_{r,k}.
\end{eqnarray*}
As before, we obtain a system of 3 linear equations on the variables
$f'_{r,k+m},\,f'_{r,k},\,f'_{r,k-m}$, whose  determinant of
coefficients is
\begin{eqnarray*}
\D_2(m)\!\!\!&=&\!\!\!-(b-b')(b-b'+1)(b-b'+2)(b+b'+2)(b^2+2bb'+5b+3b'+{b'}^2+3)m^6,
\end{eqnarray*}
and thus
\begin{eqnarray}\label{eqconddd1}
&&b'\in\{-b-2,b,b+1,b+2,\frac{1}{2}(-3-2b\pm \sqrt{-3-8b})\}.
\end{eqnarray}
Now replacing $k\in\Z$ by $k\in\frac12+\Z$ and repeating the above
arguments, we obtain
\begin{eqnarray}\label{040801}
&&b'\in\{-b-1,b,b+1,b+2,\frac{1}{2}(-1-2b\pm \sqrt{1-8b})\}\mbox{, \
\ or }\nonumber\\&& b'\in\{-b-2,b-1,b-2,b,\frac{1}{2}(-5-2b\pm
\sqrt{13+8b})\}.
\end{eqnarray}
For convenience, we introduce the following functions:
\begin{eqnarray*}
\Na_{1}(x,y)\!\!\!&=&\!\!\!
(2x+2y+3)(4+3x-3x^2-2x^3+12y\\
&&\!\!\!+4xy-2x^2y+9y^2+2xy^2+2y^3),\label{Na1}\\
\Na_{2}(x,y)\!\!\!&=&\!\!\!
18(x+y+1)(x+y+2)(a-k),\label{Na2}\\
\Na_{3}(x,y)\!\!\!&=&\!\!\!
4(-12-23x-12x^2+x^4-32y-33xy\nonumber\\
&&\!\!\!-5x^2y+2x^3y-27{y}^2-14b{b'}^2-9{y}^3-2x{y}^3-y^2).\label{Na3}
\end{eqnarray*}
For $m,\,n,\,k,\,p\in\Z$, applying
\begin{eqnarray}\label{eqactionLLG}
&&(\frac{m+n}{2}-p)\big[L_m,[L_n,G_p]\big]=(\frac{n}{2}-p)(\frac{m}{2}-n-p)[L_{m+n},G_p]
\end{eqnarray}
to $x_k$, we can as before obtain a system of 3 linear equation on
$g_{p,k-m},g_{p,k},g_{p+m}$, whose determinant of coefficients is
\begin{eqnarray*} \label{De32m1}
\!\!\!\Delta_3(m,p,k)\!\!\!&=&\!\!\!-\frac{m^6p}{4}(b-b'-1)(b-b')
%\nonumber\\&&\times
(\Na_{1}(b,b')m^2+\Na_{2}(b,b')p+\Na_{3}(b,b')p^2).
\end{eqnarray*}
Thus Lemma \ref{lemma2.1}(iii) shows $\Delta_3(m,p,k)=0$ for all
$m,p,k\in\Z$, which implies $b'\in\{b-1,b\}$, or
$\Na_{1}(b,b')=\Na_{2}(b,b')=\Na_{3}(b,b')=0, $
i.e.,
\begin{eqnarray}\label{040803} b'\in\{b-1,b\},\ \,{\rm or}\,\
\,(b,b')\in\Om=
\{(-\frac{3}{2},\,-\frac{1}{2}),\,(-1,\,0),\,(0,\,-2),
\,(\frac{1}{2},\,-\frac{3}{2})\}.
\end{eqnarray}
Similarly, using $g'_{p,m}$ instead of $g_{p,m}$, we obtain
\begin{eqnarray}\label{040804}
b'\in\{b+1,b\},\ \,{\rm or}\ \,(b,b')\in\Om'=
\{(-\frac{1}{2},\,-\frac{3}{2}),\,(0,\,-1),\,(-2,\,0),
\,(-\frac{3}{2},\,\frac{1}{2})\}.
\end{eqnarray}
Using (\ref{040801}), (\ref{040803}) and (\ref{040804}), and that
one of (\ref{Delta12=0}) and (\ref{eqconddd1}) holds, we obtain the
lemma. \QED \vs{6pt}

We shall denote the module $A_{a,b,b}$ by $A_{a,b}$ for convenience
in the following. Now we begin to prove Theorem \ref{mainth}. Noting
the fact that $L_0=G^2_0$ and the action of $L_0$ on the module
$A_{a,b}$, we claim that $g_{0,k}\neq0$ for most
$k\in\frac{1}{2}\Z$. We need to determine
$f_{r,k},\,f'_{r,k},\,g_{n,k},\,g'_{n,k}$ defined in
(\ref{modstr102}), where $n\in\Z$, $r\in\frac{1}{2}+\Z$ and
$k\in\frac{1}{2}\Z$. From the relation
$[G_p,G_n]x_k=(p-n)T_{p+n}x_k$ for
$n\in\Z,\,p\in\frac{1}{2}+\Z,\,k\in\frac{1}{2}\Z$, we only need to
determine $g_{n,k}$ and $g'_{n,k}$ for all
$n\in\Z,k\in\frac{1}{2}\Z$.

\begin{lemm}\adddot\label{Lemma-case1}
For the module $A_{a,b}$, one has
\begin{eqnarray*}
&&g_{n,k+m}=g_{n,k},\
\,\big(a-k+(2b+1)n\big)g'_{n,k+m}=\big(a-k-m+(2b+1)n\big)g'_{n,k},
\end{eqnarray*}
for $m,n\in\Z$, $k\in\frac{1}{2}\Z$.
\end{lemm}
{\bf Proof.}\ \ For  $m,\,n,\,k,\,p\in\Z$, applying
(\ref{eqactionLLG}) to $x_k$, comparing the coefficients of
$y_{k+m+n+p}$ and replacing $m,\,n,\,k$ by (i) $m,\,m,\,k-m$, (ii)
$-m,\,-m,\,k+m$, we obtain two equations. Then canceling $g_{k-m}$,
one has
\begin{eqnarray*}
&&\big(\Delta_1(b)m^4+\Delta_1(k,p)m^2+\Delta'_1(k,p)\big)(g_{p,k}-g_{p,k+m})=0,
\mbox{ \ \ where}\\
&&\Delta_1(b)=-(3+13b+18b^2+8b^3),\\
&&\Delta_1(k,p)=(3+4b)(a-k)^2+(2b^2+7b+6)(a-k)p+2(6+19b+23b^2+10b^3)p^2,\\
&&\Delta'_1(k,p)=2p\big((3(a-k)^3-2p(b+3)(a-k)^2-2bp^2
(5+4b)(a-k)+4b(b+1)p^3\big),
\end{eqnarray*}
which implies $g_{n,k+m}=g_{n,k}$, $\forall\,m,n,k\in\Z$. Similarly,
$g_{n,k+m}=g_{n,k}$, $\forall\,m,n\in\Z,k\in\frac{1}{2}+\Z$.

Replacing $x_k$ with $y_k$ in the above process, one has
\begin{eqnarray*}
&&\big(2(3+2b)p^3+4(k-a)(3+2b)p^2+(3+2b)(1+3b)pm^2+6(a-k)^2p\\
&&+(k-a)(3+4b)m^2\big)\big((a-k+2b p+p)g'_{p,k+m} -(a-k-m+2b
p+p)g'_{p,k}\big)=0,
\end{eqnarray*}
which implies
$(a-k+2bn+n)g'_{n,k+m}=(a-k-m+2bn+n)g'_{n,k},\,\forall\,m,n,k\in\Z$.

Replacing $x_k$ with $y_k$, $k\in\Z$ with $k\in\frac12+\Z$ and
repeating the above process, one has
\begin{eqnarray*}
&&\big(2(3+2b)p^3+4(k-a)(3+2b)p^2+(3+2b)(1+3b)pm^2+6(a-k)^2p\\
&&+(k-a)(3+4b)m^2\big)\big((a-k+2b p+p)g'_{p,k+m} -(a-k-m+2b
p+p)g'_{p,k}\big)=0,
\end{eqnarray*}
which implies
$(a-k+2bn+n)g'_{n,k+m}=(a-k-m+2bn+n)g'_{n,k},\,\forall\,m,n\in\Z,k\in\frac{1}{2}+\Z$.
The  lemma follows.\QED

\begin{lemm}\adddot\label{eqfcla2.6}
For the module $A_{a,b}$, one has
\begin{eqnarray*}
&&g_{n,k}=\left\{\begin{array}{ll}
\!\!\alpha_1&{\rm if}\ \,k\in\Z,\vs{6pt}\\
\!\!\alpha_2&{\rm if}\ \,k\in\frac{1}{2}+\Z,
\end{array}\right.
\ \ g'_{n,k}=\left\{\begin{array}{ll}
\!\!(a-k+2bn+n)\alpha_3&{\rm if}\ \,k\in\Z,\vs{6pt}\\
\!\!(a-k+2bn+n)\alpha_4&{\rm if}\ \,k\in\frac{1}{2}+\Z,
\end{array}\right.
\end{eqnarray*}
for $n\in\Z$ and some $\alpha_i\in\C$, $i=1,\cdots,4$.
\end{lemm}
{\bf Proof.}\ \ For $m,n,k\in\Z$, applying
$[L_m,G_n]=(\frac{m}{2}-n)G_{m+n}$ to $x_k$ and comparing the
coefficients of $y_{k+m+n}$, one has
\begin{eqnarray}\label{eqfnkfnm1}
&&(a-k-n+bm+\frac{m}{2})g_{n,k}-(a-k+bm)g_{n,m+k}=(\frac{m}{2}-n)g_{m+n,k},
\end{eqnarray}
which together with Lemma \ref{Lemma-case1}, gives
\begin{eqnarray}\label{eqfnk??m1??}
&&(\frac{m}{2}-n)(g_{m+n,k}-g_{n,k})=0.
\end{eqnarray}
From this and Lemma \ref{Lemma-case1}, we obtain
\begin{eqnarray}\label{eqfk+mfnk2}
&&g_{n,k}=\alpha_1,\ \,\forall\,\,n,k\in\Z.
\end{eqnarray}
For $m,n\in\Z,k\in\frac{1}{2}+\Z$, applying
$[L_m,G_n]=(\frac{m}{2}-n)G_{m+n}$ to $x_k$ and comparing the
coefficients of $y_{k+m+n}$ and using Lemma \ref{Lemma-case1}, we
obtain
\begin{eqnarray}\label{041001}
&&g_{n,k}=\alpha_2,\ \,\forall\,\,n\in\Z,k\in\frac{1}{2}+\Z.
\end{eqnarray}
For $m,n,k\in\Z$, applying $[L_m,G_n]=(\frac{m}{2}-n)G_{m+n}$ to
$y_k$ and comparing the coefficients of $x_{k+m+n}$, one has
\begin{eqnarray}\label{041002}
&&(a-k-n+bm)g'_{n,k}-(a-k+bm+\frac{m}{2})g'_{n,m+k}=(\frac{m}{2}-n)g'_{m+n,k},
\end{eqnarray}
which together with Lemma \ref{Lemma-case1} (taking $l=m+n$) gives
\begin{eqnarray}\label{041003}
&&(3n-l)\big((a-k+2bn+n)g'_{l,k}-(a-k+2bl+l)g'_{n,k}\big)=0.
\end{eqnarray}
From this and Lemma \ref{Lemma-case1}, we obtain
\begin{eqnarray}\label{041004}
(3n-l)(a-k+2bl+l)\big((a-k+2bn+n)g'_{l,j}-(a-j+2bl+l)g'_{n,k}\big)=0,
\end{eqnarray}
from which we obtain (by changing their indexes or choosing suitable
indexes)
\begin{eqnarray}\label{041005}
&&(a-k+2bn+n)g'_{l,j}=(a-j+2bl+l)g'_{n,k},\
\,\forall\,\,l,j,k,n\in\Z.
\end{eqnarray}
Using this, we obtain that there are exists some $\alpha_3\in\C$
such that
\begin{eqnarray}\label{041006}
&&g'_{n,k}=(a-k+2bn+n)\alpha_3,\ \,\forall\,\,k,n\in\Z.
\end{eqnarray}
For $m,n\in\Z,k\in\frac{1}{2}+\Z$, applying
$[L_m,G_n]=(\frac{m}{2}-n)G_{m+n}$ to $y_k$ and comparing the
coefficients of $x_{k+m+n}$ and using the same discussions as those
between (\ref{041003}) and (\ref{041006}), we  obtain
\begin{eqnarray*}
&&g'_{n,k}=(a-k+2bn+n)\alpha_4,\
\,\forall\,\,n\in\Z,k\in\frac{1}{2}+\Z.
\end{eqnarray*}
The lemma follows.\QED
\begin{lemm}\adddot\label{lemm2.7}
For the module $A_{a,b}$, one has, for $r\in\frac{1}{2}+\Z$,
\begin{eqnarray*}
&&f_{r,k}=\left\{\begin{array}{ll}
\!\!\frac{1}{r}\big((a-k-r)\alpha_4-(a-k+2br+r)\alpha_1\big)
&{\rm if}\ \,k\in\Z,\vs{6pt}\\
\!\!\frac{1}{r}\big((a-k-r)\alpha_3-(a-k+2br+r)\alpha_2\big)&{\rm
if}\ \,k\in\frac{1}{2}+\Z,
\end{array}\right.\\
&&f'_{r,k}=\left\{\begin{array}{ll}
\!\!\frac{1}{r}\big((a-k)\alpha_3-(a-k+2br+r)\alpha_2\big)
&{\rm if}\ \,k\in\Z,\vs{6pt}\\
\!\!\frac{1}{r}\big((a-k)\alpha_4-(a-k+2br+r)\alpha_1\big)&{\rm if}\
\,k\in\frac{1}{2}+\Z.
\end{array}\right.
\end{eqnarray*}
\end{lemm}
{\bf Proof.}\ \ For $k\in\Z,r\in\frac{1}{2}+\Z$, applying
$T_{r}=\frac{1}{r}[G_r,G_0]$ to $x_k$ and comparing the coefficients
of $x_{k+r}$, we obtain
\begin{eqnarray*}
&&f_{r,k}=\frac{1}{r}\big((a-k-r)\alpha_4-(a-k+2br+r)\alpha_1\big),\
\ \forall\,\,r\in\frac{1}{2}+\Z,k\in\Z.
\end{eqnarray*}
The other formulas also can be similarly obtained. The lemma
follows.\QED
\begin{lemm}\adddot\label{lemm2.8}
$\alpha_i=1,\ i=1,\cdots,4. $ Thus we obtain the module $A_{a,b}$
given in $(\ref{modstr4}).$
\end{lemm} {\bf Proof.}\
\ For $r,s\in\frac{1}{2}+\Z,\,k\in\Z$, applying $[T_r,
T_s]=\frac{c}{3}r\delta_{r+s,0}$ to $x_k$ and $y_k$, then comparing
the coefficients of $x_{k+r+s}$ and $y_{k+r+s}$, one has

\begin{eqnarray*}
&&((a-k-s)\alpha_4-(a-k+2bs+s)\alpha_1)((a-k-s-r)\alpha_3-(a-k-
s+2br+r)\alpha_2)\nonumber\\
&&=((a-k-r)\alpha_4-(a-k+2br+r)\alpha_1)((a-k-s-
r)\alpha_3-(a-k-r+2bs+s)\alpha_2),\\
&&((a-k-s)\alpha_3-(a-k+2bs+s)\alpha_2)((a-k-s-r)\alpha_4-(a-k-
s+2br+r)\alpha_1)\nonumber\\
&&=((a - k - r)\alpha_3 - (a - k + 2b r +r)\alpha_2) ((a-k-s -
r)\alpha_4-(a - k - r + 2b s + s)\alpha_1).
\end{eqnarray*}
Replacing $k\in\Z$ by $k\in\frac12+\Z$ in above, one has
\begin{eqnarray*}
&&((a - k)\alpha_3- (a - k + 2b s + s)\alpha_2)((a - k - s)\alpha_4-
(a - k - s + 2b
r + r)\alpha_1)\\
&&=(a - k)\alpha_4- (a - k + 2b r + r)\alpha_1)((a - k - r)\alpha_3-
(a - k - r + 2b s + s)\alpha_2), \\[4pt]
&&((a - k)\alpha_4- (a - k + 2b s
+ s)\alpha_1)((a - k - s)\alpha_3- (a - k - s + 2b
r + r)\alpha_2)\\
&&=(a - k)\alpha_3- (a - k + 2b r + r)\alpha_2)((a - k - r)\alpha_4-
(a - k - r + 2b s + s)\alpha_1).
\end{eqnarray*}
The above four equations force $\alpha_4=\alpha_1,\
\alpha_3=\alpha_2. $
For $r,p\in\frac{1}{2}+\Z,\,k\in\Z$, applying
 $[T_r,G_p]=G_{p+r}$ to $x_k$, then comparing the
coefficients of $y_{r+p+k}$, one has $(\alpha_2-\alpha_1)(1+2b)=0. $
Similarly, by allowing $k$ to be in $k\in\frac12+\Z$ and allowing
$p$ to be in $\Z$ respectively, the same arguments give
$(\alpha_2-\alpha_1)(1+b)=0 $ and $1-\alpha_1\alpha_2=0.$ Thus
 $\alpha_2=\alpha_1=\pm1$.
If necessary, rescaling $x_{r}, r\in\frac12+\Z$ and $y_k,\,k\in\Z$
by $-1$, we can suppose $\a_i=1$.\QED \vs{10pt}

\cl{\bf\S3. \ The indecomposable Harish-Chandra module
$B_{a,b}$}\setcounter{section}{3}\setcounter{theo}{0}
\setcounter{equation}{0}

\vs{4pt}

By analysis, one can see that there is possibly another type of
module of intermediate series up to isomorphism different to
$A_{a,b,b'}$, which is temporarily denoted by $B_{a,b,b'}$ with
basis $\{x_k,y_k\,|\,k\in\frac12\Z\}$ and the following module
structures
\begin{eqnarray}\label{modstr201}\begin{array}{lllllll}
&&L_mx_i=\left\{\begin{array}{ll}
(a-i+bm)x_{m+i}&{\rm if}\ \,i\in\Z,\vs{6pt}\\
\big(a-i+m(b'+\frac{1}{2})\big)x_{m+i}&{\rm if}\
\,i\in\frac{1}{2}+\Z,
\end{array}\right.\vs{6pt}\\
&&L_my_j=\left\{\begin{array}{ll}
(a-j+b'm)y_{m+j}&{\rm if}\ \,j\in\Z,\vs{6pt}\\
\big(a-j+m(b+\frac{1}{2})\big)y_{m+j}&{\rm if}\
\,j\in\frac{1}{2}+\Z,
\end{array}\right.\vs{6pt}\\
&&G_px_i=\left\{\begin{array}{cl}
y_{p+i}&{\rm if}\ \,i\in\Z,\vs{6pt}\\
-\big(a-i+2p(b'+\frac{1}{2})\big)y_{p+i}&{\rm if}\
\,i\in\frac{1}{2}+\Z,
\end{array}\right.\vs{6pt}\\
&&G_py_i=\left\{\begin{array}{cl}
x_{p+i}&{\rm if}\ \,i\in\Z,\vs{6pt}\\
-\big(a-i+2p(b+\frac{1}{2})\big)x_{p+i}&{\rm if}\
\,i\in\frac{1}{2}+\Z,
\end{array}\right.\\[18pt]
&&T_rx_k=f_{r,k}x_{k+r},\ \,\,\ T_ry_k=f'_{r,k}y_{k+r}, \ \ \
G_nx_k=g_{n,k}y_{k+n},\ \,\ G_ny_k=g'_{n,k}x_{k+n},
\end{array}
\end{eqnarray}
where $n\in\Z,\,r,p\in\frac{1}{2}+\Z,\,i,j,k\in\frac{1}{2}\Z,
\,a,\,b,\,b',\,f_{r,k},\,f'_{r,k},\,g_{n,k},\,g'_{n,k}\in{\ma C}$.
\begin{lemm}\adddot\label{lemma2.4}
$b'=b\pm\frac{1}{2}$ or $(b,b')\in\{(-\frac{3}{2},\,0),
\,(0,\,-\frac{3}{2})\}.$
\end{lemm}
{\bf Proof.}\ \ Using the similar discussions given in
(\ref{supp040701})--(\ref{040801}), one has
$b'\in(\Lambda_1\cup\Lambda_2)\cap(\Lambda_3\cup\Lambda_4)$, where
\begin{eqnarray*}
&&\Lambda_1=\{-\frac{3}{2}\pm b,b-\frac{5}{2},b-\frac{1}{2},-2-b\pm
\frac{1}{2}\sqrt{9+8b}\},\\
&&\Lambda_2=\{\pm(b+\frac{3}{2}),b+\frac{5}{2},b+\frac{1}{2},-1-b\pm
\frac{1}{2}\sqrt{-3-8b}\},\\
&&\Lambda_3=\{\pm(b+\frac{3}{2}),b\pm\frac{1}{2},-1-b\pm
\frac{1}{2}\sqrt{1-8b}\},\\
&&\Lambda_4=\{-\frac{3}{2}\pm b,b\pm\frac{1}{2},-2-b\pm
\frac{1}{2}\sqrt{13+8b}\}.
\end{eqnarray*}
Using the similar discussions given in
(\ref{eqactionLLG})--(\ref{040804}), we obtain $b'=b\pm\frac{1}{2}$
or $(b,b')\in\Lambda\cap\Lambda'$, where
$\Lambda=\{(-\frac{3}{2},0),(-1,\frac{1}{2}),(\frac{1}{2},-1),
(0,-\frac{3}{2})\}$, $\Lambda'=\{(-2,-\frac{1}{2}),
(-\frac{3}{2},0),(-\frac{1}{2},-2),(0,-\frac{3}{2})\}$. Then the
lemma follows.\QED
\begin{rema}\rm
(1) The modules $B_{a,0,-\frac{3}{2}}$ can be regarded as
$B_{a,-\frac{3}{2},0}$ by exchanging $x_i$'s and $y_i$'s. We shall
consider the module $B_{a,0,-\frac{3}{2}}$.
\\
(2) The modules $B_{a,b,b+\frac{1}{2}}$ can be regarded as
$B_{a,b',b'-\frac{1}{2}}$ with $b'=b+\frac12$ by exchanging $x_i$'s
and $y_i$'s. We shall consider the module $B_{a,b,b-\frac{1}{2}}$,
 and denote it as $B_{a,b}$ for convenience.
\end{rema}

\begin{lemm}\adddot\label{Lemma-cas2}
For the module $B_{a,b}$, one has\\
$(1)$\ \ if $b'=b-\frac{1}{2}$, then
\begin{eqnarray*}
&&\big(a-k+2bn\big)g_{n,k+m}=\big(a-k-m+2bn\big)g_{n,k},\,\
g_{n,k'+m}=g_{n,k'},\\
&&g'_{n,k+m}=g'_{n,k},\
\,\big(a-k'+(2b+1)n\big)g'_{n,k'+m}=\big(a-k'-m+(2b+1)n\big)g'_{n,k'},
\end{eqnarray*}
$(2)$\ \ if $b=0,\,b'=-\frac{3}{2}$, then
\begin{eqnarray*}
&&(a - k')g_{n, k'}=(a-k'-m)g_{n, k'+ m},\\
&&(a - k - n)g'_{n, k}=(a-k-m-n)g'_{n, k + m},\\
&&(a - k - m)(a - k - m - 2 n)g_{n, k}=(a - k)(a - k - 2 n)g_{n, k +
m},\\
&&(a-k'-m-n)(a-k'-m+n)g'_{n,k'}=(a-k'-n)(a-k'+n)g_{n,k'+m},
\end{eqnarray*}
for $k\in\Z,\,k'\in\frac{1}{2}+\Z$.
\end{lemm}
{\bf Proof.}\ \ (1) For $m,\,n,\,k,\,p\in\Z$, applying
(\ref{eqactionLLG}) to $x_k$, comparing the coefficients of
$y_{k+m+n+p}$ and replacing $m,\,n,\,k$ by (i) $m,\,m,\,k-m$, (ii)
$-m,\,-m,\,k+m$, we obtain two equations. Canceling $g_{k-m}$ and
$b'=b-\frac{1}{2}$, one has
\begin{eqnarray*}
&&\big(4(1+b)p^3+8(k-a)(1+b)p^2+(1+b)(6b-1)pm^2+6(a-k)^2p\\
&&+(k-a)(1+4b)m^2\big)\big((a-k+2b p)g'_{p,k+m} -(a-k-m+2b
p)g'_{p,k}\big)=0,
\end{eqnarray*}
which implies
$(a-k+2bn)g'_{n,k+m}=(a-k-m+2bn)g'_{n,k},\,\forall\,m,n,k\in\Z$.
Replacing $x_k$ by $y_k$ in above, we obtain
\begin{eqnarray*}
&&\!\!\!\!\!\!\big(\Delta_2(b)m^4+\Delta_2(k,p)m^2+\Delta'_2(k,p)\big)(g'_{p,k}-g'_{p,k+m})=0,
\mbox{ \ \ where}\\
&&\!\!\!\!\!\!\Delta_2(b)=-(b+6b^2+8b^3),\\
&&\!\!\!\!\!\!\Delta_2(k,p)=(1+4b)(a-k)^2+(2b^2+5b+3)(k-a)p+(2+7b+16b^2+20b^3)p^2,\\
&&\!\!\!\!\!\!\Delta'_2(k,p)=2p\big(3(a-k)^3-p(5+2b)(a-k)^2+(1-2b)
(3+4b)(a-k)p^2+2b(2b-1)p^3\big).
\end{eqnarray*}
Hence $g'_{n,k+m}=g'_{n,k},\,\forall\,m,n,k\in\Z$. Similarly,
$g_{n,k'+m}=g_{n,k'},\,\forall\,m,n\in\Z,k'\in\frac{1}{2}+\Z$.

As above, for $m,\,n,\,p\in\Z,\,k'\in\frac{1}{2}+\Z$, applying
(\ref{eqactionLLG}) to $y_{k'}$, one can obtain
\begin{eqnarray*}
&&\big(2(3+2b)p^3+4(k'-a)(3+2b)p^2+(3+2b)(1+3b)pm^2+6(a-k')^2p\\
&&+(k'-a)(3+4b)m^2\big)\big((a-k'+2b p+p)g'_{p,k'+m} -(a-k'-m+2b
p+p)g'_{p,k'}\big)=0,
\end{eqnarray*}
which implies
$(a-k'+2bn+n)g'_{n,k'+m}=(a-k'-m+2bn+n)g'_{n,k'},\,\forall\,m,n\in\Z,k'\in\frac{1}{2}+\Z$.

(2) can be obtained similarly. \QED

\begin{lemm}\adddot\label{lemm3.4}
For the module $B_{a,b}$, one has\\
(1)\ \,if $b'=b-\frac{1}{2}$, then
\begin{eqnarray*}
&&\!\!\!\!\!\!\!g_{n,k}=\left\{\begin{array}{cl}
\!\!(a-k+2bn)\beta_1&{\rm if}\ \,k\in\Z,\vs{6pt}\\
\!\!\beta_2&{\rm if}\ \,k\in\frac{1}{2}+\Z,
\end{array}\right.
g'_{n,k}=\left\{\begin{array}{cl}
\!\!\beta_3&{\rm if}\ \,k\in\Z,\vs{6pt}\\
\!\!(a-k+2bn+n)\beta_4&{\rm if}\ \,k\in\frac{1}{2}+\Z,
\end{array}\right.
\end{eqnarray*}
(2)\ \ if $b=0,\,b'=-\frac{3}{2}$, then
\begin{eqnarray*}
&&\!\!\!\!\!\!\!g_{n,k}=\left\{\begin{array}{cl}
\!\!(a - k) (a - k - 2 n)\mu_1&{\rm if}\ \,k\in\Z,\vs{6pt}\\
\!\!\frac{1}{a-k}\mu_2&{\rm if}\
\,a\notin\frac{1}{2}+\Z,\,\,k\in\frac{1}{2}+\Z,
\end{array}\right.\\
&&\!\!\!\!\!\!\!g'_{n,k}=\left\{\begin{array}{cl}
\!\!\frac{1}{a-k-n}\mu_3&{\rm if}\ \,a\notin\Z,\,\,k\in\Z,\vs{6pt}\\
\!\!(a-k-n)(a-k+n)\mu_4&{\rm if}\ \,k\in\frac{1}{2}+\Z,
\end{array}\right.
\end{eqnarray*}
for $n\in\Z$ and some $\beta_i,\mu_i\in\C$, $i=1,\cdots,4$.
\end{lemm}
{\bf Proof.}\ \ (1) For $m,n,k\in\Z$, applying
$[L_m,G_n]=(\frac{m}{2}-n)G_{m+n}$ to $x_k$ and comparing the
coefficients of $y_{k+m+n}$, one has
\begin{eqnarray}\label{080410011}
&&(a-k-n+bm-\frac{m}{2})g_{n,k}-(a-k+bm)g_{n,m+k}=(\frac{m}{2}-n)g_{m+n,k},
\end{eqnarray}
which together with Lemma \ref{Lemma-cas2}, gives
\begin{eqnarray}\label{080410012}
&&(3n-j)\big((a-k+2bj)g_{n,k}-(a-k+2bn)g_{j,k}\big)=0.
\end{eqnarray}
From this and  Lemma \ref{Lemma-cas2}, we obtain
$(a-k+2bj)\big((a-i+2bj)g_{n,k}-(a-k+2bn)g_{j,i}\big)=0$ for
$n,i,j,k\in\Z$, which implies
\begin{eqnarray*}
&&g_{n,k}=(a-k+2bn)\beta_1,\ \,\forall\,\,n,k\in\Z.
\end{eqnarray*}
For $m,n\in\Z,k\in\frac{1}{2}+\Z$, applying
$[L_m,G_n]=(\frac{m}{2}-n)G_{m+n}$ to $x_k$ and comparing the
coefficients of $y_{k+m+n}$, we  obtain that (\ref{eqfnkfnm1}) and
(\ref{eqfnk??m1??}) also hold for $k\in\frac12+\Z$. From this and
Lemma \ref{Lemma-cas2}, we obtain
\begin{eqnarray}\label{08041001}
&&g_{n,k}=\beta_2,\ \,\forall\,\,n\in\Z,k\in\frac{1}{2}+\Z.
\end{eqnarray}
For $m,n,k\in\Z$, applying $[L_m,G_n]=(\frac{m}{2}-n)G_{m+n}$ to
$y_k$ and comparing the coefficients of $x_{k+m+n}$, one has
\begin{eqnarray}\label{08041002}
&&(a-k-n+bm)g'_{n,k}-(a-k+bm-\frac{m}{2})g'_{n,m+k}=(\frac{m}{2}-n)g'_{m+n,k},
\end{eqnarray}
which together with Lemma \ref{Lemma-cas2}, gives
\begin{eqnarray}\label{08041003}
&&(\frac{m}{2}-n)(g'_{m+n,k}-g'_{n,k})=0.
\end{eqnarray}
From this and  Lemma \ref{Lemma-cas2}, we obtain
\begin{eqnarray}\label{08041004}
&&g'_{n,k}=\beta_3,\ \,\forall\,\,k,n\in\Z.
\end{eqnarray}
For $m,n\in\Z,k\in\frac{1}{2}+\Z$, applying
$[L_m,G_n]=(\frac{m}{2}-n)G_{m+n}$ to $y_k$ and comparing the
coefficients of $x_{k+m+n}$, one has
$$(a-k-n+bm)g'_{n,k}-(a-k+bm+\frac{m}{2})g'_{n,m+k}=(\frac{m}{2}-n)g'_{m+n,k},
$$ which implies $g'_{n,k}=(a-k+2bn+n)\beta_4,\,\forall\,n\in\Z,k\in\frac{1}{2}+\Z. $ Then $(1)$ follows.

(2) is obtained similarly.\QED
\begin{rema}\adddot\label{rema3.5}\rm
If $b=0,\,b'=-\frac{3}{2}$, one can deduce from the above two lemmas
that
\begin{eqnarray}\label{0522a01}
&&g_{n,k}=0\ \ \mbox{if}\ \ a,k\in\frac{1}{2}+\Z,\ \ \mbox{and}\ \
g'_{n,k}=0\ \ \mbox{if}\ \ a,k\in\Z.
\end{eqnarray}
However, from $G_n^2=2L_{2n}$ and $G_p^2=-2L_{2p}$ for
$n\in\Z,p\in\frac{1}{2}+\Z$, we see that (\ref{0522a01}) cannot
happen. Hence in the module $B_{a,0,-\frac{3}{2}}$, we always
suppose $a\notin\frac{1}{2}\Z$.
\end{rema}
Using the similar techniques to those used in Lemma \ref{lemm2.8},
one can obtain the following lemma.
\begin{lemm}\adddot\label{lemm3.6}
For the module $B_{a,b}$, one has\\
(1)\ \,if $b'=b-\frac{1}{2}$, then
\begin{eqnarray*}
&&\!\!\!\!\!\!\!\!f_{r,k}=\left\{\begin{array}{cl}
\!\!\frac{1}{r}\big((a-k)\beta_1+(a-k-r)\beta_4\big)
&{\rm if}\ \,k\in\Z,\vs{6pt}\\
\!\!-\frac{1}{r}\big((a-k+2br+r)\beta_2+(a-k+2br)\beta_3\big)&{\rm
if}\ \,k\in\frac{1}{2}+\Z,
\end{array}\right.\\
&&\!\!\!\!\!\!\!\!f'_{r,k}=\left\{\begin{array}{cl}
\!\!\frac{1}{r}(\beta_2+\beta_3)
&{\rm if}\ \,k\in\Z,\vs{6pt}\\
\!\!-\frac{1}{r}\big((a-k)(a-k+2br)\beta_4+(a-k-r)(a-k+2br+r)\beta_1\big)&{\rm
if}\ \,k\in\frac{1}{2}+\Z,
\end{array}\right.
\end{eqnarray*}
(2)\ \ if $a\notin\frac{1}{2}\Z$, $b=0,\,b'=-\frac{3}{2}$, then
\begin{eqnarray*}
&&\!\!\!\!\!\!\!\!f_{r,k}=\left\{\begin{array}{cl}
\!\!\frac{1}{r}(a-k)^2(\mu_1+\mu_4)
&{\rm if}\ \,k\in\Z,\vs{6pt}\\
\!\!-\frac{1}{r(a-k)}\big((a-k-2r)\mu_3+(a-k+r)\mu_2\big)&{\rm if}\
\,k\in\frac{1}{2}+\Z,
\end{array}\right.\\
&&\!\!\!\!\!\!\!\!f'_{r,k}=\left\{\begin{array}{cl}
\!\!\frac{1}{r(a-k)}(\mu_2+\mu_3)
&{\rm if}\ \,k\in\Z,\vs{6pt}\\
\!\!-\frac{1}{r}\big((a-k+r)(a-k-r)^2\mu_1+(a-k)^2(a-k-2r)\mu_4\big)&{\rm
if}\ \,k\in\frac{1}{2}+\Z,
\end{array}\right.
\end{eqnarray*}
for $r\in\frac{1}{2}+\Z$.
\end{lemm}

Similar to the proof of Lemma \ref{lemm2.8}, we obtain
\begin{lemm}\adddot\label{lemm11}
$(1)\ \,\beta_1=-\beta_2= \beta_3=-\beta_4=1;$ \ \ $(2)\ \,\mu_i=0,\
i=1,\cdots,4.$
\end{lemm}
\begin{rema}\adddot\label{rema2.14}\rm
The above two lemmas show that $g_{n,k}= g'_{n,k}=0$ for $
k\in\frac{1}{2}\Z$ in the module $B_{a,0,-\frac{3}{2}}$, thus such a
module  does not exist for any $a\in\C$. Hence we obtain the module
$B_{a,b}$ defined in (\ref{modstr2}).
\end{rema}

\cl{\bf\S4. \ The deformations of the modules $A_{a,b}$ and
$B_{a,b}$}\setcounter{section}{4}\setcounter{theo}{0}
\setcounter{equation}{0}

\vs{6pt}

The proof of the main theorem will be completed by the
following.\vs{-6pt}
\begin{lemm}\adddot\label{lemm13A}
$(1)$ The module $A_{a,b}$ has two types deformations up to
isomorphism, denoted by $A_1(\a)$ and $A_2(\a)$, whose module
structures are respectively determined by $(\ref{modstr3})$ and
$(\ref{modstr4}).$\\
$(2)$ The module $B_{a,b}$ has only two types of deformations,
denoted by $B_1(\a)$ and $B_2(\a)$, whose module structures are
determined by $(\ref{modstr66})$ and $(\ref{modstr68})$
respectively.
\end{lemm}
{\bf Proof.}\ \ (1) Note that a module has a nontrivial deformation
only if it is reducible. One can see that there are two possible
cases in which the module $A_{a,b}$ has deformation.

\noindent{\bf Case 1.}\quad$a\in\Z,\,\,b=-1$.

In this case, by shifting the index $k$, we can suppose
$a=0,\,b=-1$. Then
$\{x_k\,|\,k\in(\frac12+\Z)\cup\Z^*\}\cup\{y_j\,|\,j\in\frac12\Z\}$
span the only proper submodule. Thus the possible deformations are
the actions of $T_r,G_q,L_n$ on $x_0$ for all
$n\in\Z,r\in\frac{1}{2}+\Z,q\in\frac{1}{2}\Z$. First, one has
\begin{eqnarray*}
&&L_nx_j=-(j+n)x_{n+j},\ \,\ L_ny_k=-(k+\frac{n}{2})y_{n+k},\\
&&T_rx_j=0,\ T_ry_k=y_{k+r},\ \,G_qx_j=y_{q+j},\ \,
G_qy_k=(-1)^{2q+1}(k+q)x_{q+k},\label{0418a05}
\end{eqnarray*}
for
$n\in\Z,r\in\frac{1}{2}+\Z,j\in(\frac12+\Z)\cup\Z^*,q,k\in\frac{1}{2}\Z$.
Meanwhile, one can suppose
\begin{equation*}
L_n x_0=e_n x_n,\ \ T_r x_0=f_r x_r,\ \ G_n x_0=h_n x_n,\ \ G_p
x_0=g_p x_p,
\end{equation*}
for $n\in\Z,p,r\in\frac{1}{2}+\Z$ and some $e_n,f_n,g_p,h_n\in\C$.
For $n\in\Z^*$, applying $[L_{n},L_1]=(n-1)L_{n+1}$ to $x_{0}$ and
comparing the coefficients of $x_{n+1}$, we obtain
\begin{eqnarray}\label{08041702}
&&(n+1)(e_n-e_1)=(n-1)e_{n+1},\ \ \forall\,\,n\in\Z^*,
\end{eqnarray}
which implies $e_0=0$.  Using induction on $n$ in (\ref{08041702}),
one has
\begin{eqnarray}
e_n=\left\{\begin{array}{cl}
\frac{n^2-n}{2}e_{-1}+\frac{n^2+n}{2}e_1&{\rm if}\ \,n\leq-1,\vs{6pt}\\
\frac{n^2-n}{2}e_{2}+(2n-n^2)e_1&{\rm if}\ \,n\geq0.
\end{array}\right.\label{08041703}
\end{eqnarray}
Applying $[L_{2},L_{-1}]=3L_{1}$ to $x_{0}$ and comparing the
coefficients of $x_{1}$, we obtain $e_{-1}=e_2-3e_{1}, $
using which,
one can rewrite (\ref{08041703}) as
\begin{eqnarray}\label{08041704}
e_n=-n(\a'n+\a),\ \ \forall\,\,n\in\Z,
\end{eqnarray}
where $\a=\frac{e_2}{2}-2e_1,\ \,\a'=e_1-\frac{e_2}{2}=\a-e_1$. For
$n\in\Z^*,\,p\in\frac{1}{2}+\Z$, applying
$[L_{n},G_p]=(\frac{n}{2}-p)G_{n+p}$ to $x_{0}$ and comparing the
coefficients of $y_{n+p}$, we obtain
\begin{eqnarray}\label{0422m001}
&&(p+\frac{n}{2})g_p-n(\a'n+\alpha)=(p-\frac{n}{2})g_{n+p}.
\end{eqnarray}
Letting $n=2p$ in (\ref{0422m001}), we obtain $g_p=2p\a'+\alpha,$
and also $g_{n+p}=2(n+p)\a'+\alpha.$ Taking them back to
(\ref{0422m001}), we know that $g_p=2p\a'+\alpha$ is indeed a
solution of the equation (\ref{0422m001}).

For $m,n\in\Z^*$, applying $[L_{n},G_m]=(\frac{n}{2}-m)G_{n+m}$ to
$x_{0}$ and repeating the above process, we  obtain
$h_m=2m\a'+\alpha$ for $m\in\Z^*$. For $n\in\Z^*$, applying
$[L_{n},G_0]=\frac{n}{2}G_{n}$ to $x_{0}$ and using the above
results obtained, one can deduce $h_0=\a$. Hence
\begin{eqnarray}\label{0422m002}
g_q=2q\a'+\alpha,\ \ \forall\,\,q\in\frac{1}{2}\Z.
\end{eqnarray}
For $r\in\frac{1}{2}+\Z$, according to
$T_rx_0=\frac{1}{r}[G_{r},G_0]x_0$ and using (\ref{0422m002}), one
can deduce $f_r=-2r\a'$ for $r\in\frac{1}{2}+\Z$. Then we get a
deformation of the module $A_{0,-1}$, denoted by $A_1(\a)$,  given
in (\ref{modstr3}).

\noindent{\bf Case 2.}\quad $a\in\Z,\,\,b=-\frac{1}{2}$.

In this case, by shifting the index $k$, one can suppose
$a=0,\,b=-\frac12$. Then $\C y_0$ span the only proper submodule.
Thus the possible deformations are the actions of $T_r,G_q,L_n$ on
$x_{-r},y_{-q},x_{-n}$ for all
$n\in\Z,r\in\frac{1}{2}+\Z,q\in\frac{1}{2}\Z$. Then one has
\begin{eqnarray*}
&&L_nx_{k'}=-(k'+\frac{n}{2})x_{n+k'},\ \ \ \ \,\,L_ny_i=-iy_{n+i}\\
&&T_rx_{k'}=-x_{k'+r},\ \ T_ry_k=0,\ \
G_qy_{k'}=(-1)^{2q+1}k'x_{q+k'},\ \ \ G_qx_j=y_{j+q},
\end{eqnarray*}
for $n\in\Z,r\in\frac{1}{2}+\Z,q,k',i,j,k\in\frac{1}{2}\Z$ and
$k\neq-r,i\neq-n,j\neq-q$. And we suppose
\begin{equation*}
L_n y_{-n}=e_n y_0,\ \ T_r y_{-r}=f_r y_0,\ \ G_q x_{-q}=g_q y_0,
\end{equation*}
for $n\in\Z,q,r\in\frac{1}{2}+\Z$ and some $e_n,f_n,g_q\in\C$. For
$n\in\Z^*$, applying $[L_{n},L_1]=(n-1)L_{n+1}$ to $y_{-n-1}$ and
comparing the coefficients of $y_{0}$, we obtain
\begin{eqnarray}\label{0423001}
&&(n+1)(e_n-e_1)=(n-1)e_{n+1},\ \ \forall\,\,n\in\Z^*,
\end{eqnarray}
which implies $e_0=0$.  Using induction on $n$ in (\ref{0423001}),
one also can get $e_n=n(\a'n+\a),\, \forall\,\,n\in\Z, $ where
$\a=2e_1-\frac{e_2}{2},\ \,\a'=\frac{e_2}{2}-e_1=e_1-\a$. For
$n\in\Z^*,\,p\in\frac{1}{2}+\Z$, applying
$[L_{n},G_p]=(\frac{n}{2}-p)G_{n+p}$ to $x_{-n-p}$ and comparing the
coefficients of $y_{0}$, we obtain
\begin{eqnarray}\label{0423004}
&&(p+\frac{n}{2})g_p-n(\a'n+\alpha)=(p-\frac{n}{2})g_{n+p}.
\end{eqnarray}
Letting $n=2p$ in (\ref{0423004}), we obtain $g_p=2p\a'+\alpha,$ and
also $g_{n+p}=2(n+p)\a'+\alpha.$ Taking them back to
(\ref{0423004}), we know that $g_p=2p\a'+\alpha$ is indeed the
solution of the equation (\ref{0423004}). For $m,n\in\Z^*$, applying
$[L_{n},G_m]=(\frac{n}{2}-m)G_{n+m}$ to $x_{-n-m}$ and repeating the
above process, we also obtain $g_m=2m\a'+\alpha$ for $m\in\Z^*$. For
$n\in\Z^*$, applying $[L_{n},G_0]=\frac{n}{2}G_{n}$ to $x_{-n}$ and
using the above results obtained, one also can deduce $g_0=\a$.

For $r\in\frac{1}{2}+\Z$, according to
$T_ry_{-r}=\frac{1}{r}[G_{r},G_0]y_{-r}$, one can deduce $f_r=2\a'r$
for $r\in\frac{1}{2}+\Z$. Then we get a deformation of the module
$A_{0,-\frac{1}{2}}$, denoted by ${A}_2(\a)$,  given in
(\ref{modstr4}).

The proof of (2) is similar. This proves the lemma and Theorem
\ref{mainth}.\QED

\end{document}